\documentclass{article}
\usepackage{graphicx}
\usepackage{amsmath,amsfonts,amsthm}
\usepackage{amssymb,latexsym}
\usepackage{fixmath}
\usepackage{mathrsfs,amsbsy}
\usepackage{dsfont}
\usepackage{enumerate}
\usepackage{latexcad}
\usepackage{ulem}
\usepackage{eucal}

\usepackage{hyperref}
\hypersetup{colorlinks}

\usepackage{colortbl}
\usepackage{multirow}

\def\wtd{\widetilde}
\def\what{\widehat}

\def\cK{{\cal K}}
\def\cN{{\cal N}}
\def\cR{{\cal R}}

\def\bb{\pmb{b}}

\def\be{\pmb{e}}

\def\bg{\pmb{g}}

\def\bq{\pmb{q}}

\def\bx{\pmb{x}}

\DeclareMathOperator{\diag}{diag}
\DeclareMathOperator{\eig}{eig}
\DeclareMathOperator{\ext}{ext}

\DeclareMathOperator{\SOCLCP}{SOCLCP}
\DeclareMathOperator{\SOL}{SOL}
\DeclareMathOperator{\rank}{rank}
\DeclareMathOperator{\rel}{rel}
\DeclareMathOperator{\fl}{fl}
\def\macheps{\mathfrak{u}}

\def\BN{\textsf{BN}}
\def\LvA{\textsf{LCPvA}}

\DeclareMathAlphabet{\mathpzc}{OT1}{pzc}{m}{it}

\newtheorem{proposition}{Proposition}[section]
\newtheorem{theorem}{Theorem}[section]
\newtheorem{lemma}{Lemma}[section]
\newtheorem{corollary}{Corollary}[section]

\theoremstyle{definition}

\newtheorem{remark}{Remark}[section]

\numberwithin{equation}{section}
\numberwithin{figure}{section}
\numberwithin{table}{section}

\usepackage{algorithm,algorithmic}
\numberwithin{algorithm}{section}

\def\coneK{{\mathbb K}}

\allowdisplaybreaks

\def\Red{\textcolor{red}}

\usepackage{subfigure}

\def\RKSM{\textsf{RKSM}}

\def\tdM{\tilde{M}}
\def\tdJ{\tilde{J}}
\def\tdq{\tilde{q}}
\def\tdf{\tilde{f}}

\DeclareMathOperator{\T}{T}

\DeclareMathOperator{\orth}{orth}

\def\bbR{\mathbb{R}}

\def\orth{{\rm orth}}

\title{Solving Symmetric and Positive Definite Second-Order Cone Linear Complementarity Problem
by A Rational Krylov Subspace Method%
}

\author{Yiding Lin\thanks{School of Economic  Mathematics, Southwestern University of Finance and Economics, 555 Liutai Road, Chengdu 611130, China ({\tt Yiding.Lin@gmail.com}). }
\and Xiang Wang \thanks{Department of Mathematics, Nanchang University, 999 Xuefu Road, Nanchang 330031, China
              ({\tt wangxiang49@ncu.edu.cn})}
\and Lei-Hong Zhang  \thanks{
	{Corresponding author.} School of Mathematical Sciences and Institute of Computational Science, Soochow University, Suzhou 215006, Jiangsu, China ({\tt longzlh@suda.edu.cn}).
             The work of this author  was supported in part by the National Natural Science Foundation of China
             NSFC-11671246 and NSFC-12071332.}
}

\begin{document}
\bibliographystyle{siam}
\maketitle

\begin{abstract}
The second-order cone linear complementarity problem (SOCLCP) is a generalization of the classical linear complementarity problem. It has been known that SOCLCP, with the
globally uniquely solvable  property, is essentially equivalent to a zero-finding problem in which the associated
function bears much similarity to the transfer function in model reduction
[{\em SIAM J. Sci. Comput.}, 37 (2015), pp.~A2046--A2075]. In this paper,
we propose a new rational Krylov subspace method to solve the zero-finding problem for the symmetric and positive definite SOCLCP.
The algorithm consists of two stages: first, it relies on an extended Krylov subspace to obtain rough approximations of the zero root, and then applies multiple-pole rational Krylov subspace projections iteratively to acquire an accurate solution.
Numerical evaluations on various types of SOCLCP examples demonstrate its efficiency and robustness.
\end{abstract}

{\small

}

%

\pagestyle{myheadings}
\thispagestyle{plain}

\markboth{Yiding Lin, Xiang Wang,Lei-Hong Zhang,Ren-Cang Li }
{Solving 
SPD SOCLCP	by A Rational Krylov Subspace Method}
\section{Introduction}
For a given symmetric and positive definite $M\in\bbR^{n\times n}$ and  a vector $\bq\in \bbR^{n}$, in this paper, we are concerned with the following {\em second-order cone linear complementarity problem\/}:
\begin{eqnarray}\label{eq:SOCLCP}
\SOCLCP(\coneK^n,M,\bq): \quad {\rm find}~~\bx\in \coneK^n \mbox{ such that}\quad \bq+M\bx\in\coneK^n\quad \mbox{and}\quad\bx^{\T}(\bq+M\bx)=0,
\end{eqnarray}
where $\coneK^n:=\Big\{[x_1,\bx_2^{\T}]^{\T}\in\bbR\times \bbR^{n-1}\,:\,\|\bx_2\|_2\leq x_1\Big\}$
is the so-called   {\em second-order cone\/} or the {\em the Lorentz cone\/}.
The set of solutions of $\SOCLCP(\coneK^n,M,\bq)$  is denoted by
$\SOL(M,\coneK^n,\bq)$.

The $\SOCLCP(\coneK^n,M,\bq)$ is a generalization of the classical linear complementarity problem  $\SOCLCP(\bbR_+^n,M,\bq)$ \cite{cops:1992} from  the nonnegative cone $\bbR_+^n=\{\bx\in \bbR^n|\bx\ge 0\}$ to $\coneK^n.$ Solvability and many crucial properties of the set of solutions 
of $\SOCLCP(\coneK^n,M,\bq)$ have been established (see e.g., \cite{chpa:2008,chts:2005,chss:2003,fult:2002,hayf:2005,hayf:2005b,yayu:2013,yazs:2017,zhys:2015}), which provide foundations for many numerical methods.  For example,  the notion of {\em globally uniquely solvable\/} (GUS) property for $M$ which states that the linear complementarity problem  has a unique solution for any given $\bq\in\bbR^n$, plays an important role in both $\SOCLCP(\bbR_+^n,M,\bq)$ \cite{cops:1992} and $\SOCLCP(\coneK^n,M,\bq)$ \cite{yayu:2013}. In particular, for $\SOCLCP(\coneK^n,M,\bq)$, Yang and Yuan \cite{yayu:2013}  proposed an important algebraic characterization of the GUS property. This provides new perspective on $\SOCLCP(\coneK^n,M,\bq)$  and  also leads to several efficient numerical methods (see e.g., \cite{yayu:2013,zhya:2013a,zhya:2014a,zhys:2015}). A recent work 
\cite{walz:2019} is an factorization-based numerical algorithm that is efficient for   small- to medium-size problems. Specifically, the method is based on the full eigen-decomposition of matrix pencil $M-\lambda J$, where
\begin{equation}\label{eq:Jn}
J=\diag(1,-1,\ldots,-1),
\end{equation}
and can be regarded as a direct method, as opposed to
iterative methods previously,  because computing the full eigen-decomposition, although iterative in nature, is mature enough to be widely considered as direct in applications \cite{demm:1997,govl:2013}. Nonetheless, for large scale and sparse problems,  such an eigen-decomposition-based method is still very expensive or inapplicable.

The focus of this paper is on the efficient methods for large-scale  $\SOCLCP(\coneK^n,M,\bq)$ with $M$ symmetric and positive definite. Our approach follows the same numerical framework  of the Krylov subspace method in \cite{zhys:2015}, which generally contains the following three steps:
\begin{itemize}
\item[{\it step 1.}] transform $\SOCLCP(\coneK^n,M,\bq)$ to a zero-finding problem \cite{zhys:2015},
\item[{\it step 2.}] project the zero-finding problem
to a much smaller scale problem by certain rational Krylov subspaces, 
\item[{\it step 3.}] solve the projected problem by an efficient zero-finding solver, e.g., the direct method of  \cite{walz:2019}.
\end{itemize}
Our main contribution is on the second step in which we provide a more effective rational Krylov subspace to improve the performance of the algorithmic framework of \cite{zhys:2015}.  Even though the method of \cite{zhys:2015} is applicable for more general $\SOCLCP(\coneK^n,M,\bq)$ where $M$ is not necessarily symmetric and positive definite, we will demonstrate that our new Krylov subspace method performs generally better than that of \cite{zhys:2015} whenever $M$ is symmetric and positive definite. 

{
This paper is organized as follows. In Section \ref{section:preliminaries}, we summarize some preliminary theoretical results of  $\SOCLCP(\coneK^n,M,\bq)$, most of which are from \cite{zhys:2015}.
In Section \ref{section:RKSM}, we propose our rational Krylov subspace projection method. The numerical approach for  the reduced problem in {\it step 3} is discussed in Subsection \ref{ssec:th(s)=0}.
Our numerical evaluation of the new method is carried out in Section \ref{section:numerical_experiments}, where we report our numerical experiments by comparing it with other algorithms. We draw our final remarks  in Section \ref{section:conclusion}.
}
\vskip 2mm

{\bf Notation}:
Throughout this paper, all
vectors are column vectors and are typeset in bold lower case letters. The identity matrix in $\bbR^{n\times n}$ will be denoted by $I_n\equiv [\be_1,\be_2,\dots,\be_n]$, where $\be_i$ is its $i$ column.
For $A\in \bbR^{n\times m}$, $A^{\T}$ denotes its transpose. The labels
$\cR(A)$ and $\cN(A)$ denote the range and kernel of $A$, respectively. Thus, $\cR(A)^\perp=\cN(A^{\T}),$ where $\cR(A)^\perp$ denotes the orthogonal complement of $\cR(A)$.
If $A$ is square, then we use $\eig(A)=\{\lambda_i(A)\}_{i=1}^n$ to represent the set of eigenvalues, and use $A\succ 0$ to indicate $A$ is symmetric and positive definite.
 Also, both $\orth(A)$ and $\orth(\cR(A))$ represent the orthonormal basis matrix of $\cR(A)$.  For convenience,  we shall adopt MATLAB-like format:
 $\bx_{(i)}$ is  the $i$th element of $\bx$ and $A_{(i,j)}$ is the $(i,j)$th entry of $A$, where $(i:j)$ stands for the set of integers from $i$ to $j$
inclusive, and $A_{(k:\ell,i:j)}$ is the  submatrix of $A$ that consists of intersections from row $k$ to
row $\ell$ and column $i$ to column $j$.
The $\ell$th Krylov subspace generated by $A\in\bbR^{n\times n}$ on $\bx\in \bbR^{n}$ is defined as
\[
\cK_{\ell}(A,\bx)={\rm span}(\bx,A\bx,\dots,A^{\ell-1}\bx).
\]
When $A$ is also invertible, the $(\ell,k)$th extended Krylov subspace \cite{knsi:2011,simo:2007}
of $A$ on $\bx$ is defined as
$$
\cK_{\ell,k}^{\ext}(A,\bx)={\rm span}(\bx,A\bx,\dots,A^{\ell-1}\bx; A^{-1}\bx,A\bx,A^{-2}\bx,\dots,A^{k-1}\bx,A^{k}\bx).
$$

%
%

\section{Preliminaries}\label{section:preliminaries}
In this section, we first briefly review preliminary results on $\SOCLCP(\coneK^n,M,\bq)$ \eqref{eq:SOCLCP}.  We begin with  a characterization of the solution which describes three mutually exclusive cases for the solution of
$\SOCLCP(\coneK^n,M,\bq)$.

\begin{theorem}[\cite{zhya:2014a}]\label{thm:3m-e-c}
There are three mutually exclusive cases
for the solution $\bx\in \SOL(\coneK^n,M,\bq)$, namely:
\begin{itemize}
\item[{\rm (C1)}] $\bq\in \coneK^n$ (which implies that $\bx=\mathbf{0}$ is the solution);
\item[{\rm (C2)}] $\SOL(\coneK^n,M,\bq)\supseteq\{\bx\in \coneK^n\,:\,M\bx+\bq=\mathbf{0}\}\neq \varnothing$;
\item[{\rm (C3)}] there exists $s_*>0$ such that $M\bx+\bq=s_*J_n\bx\in \partial(\coneK^n)$, where
$\partial(\coneK^n)=\Big\{[x_1,\bx_2^{\T}]^{\T}\in\bbR\times \bbR^{n-1}\,:\,\|\bx_2\|_2= x_1\Big\}$ denotes the boundary of $\coneK^n$, and $J_n$ is given in \eqref{eq:Jn}.
\end{itemize}
\end{theorem}

Note that the first case (C1) is trivially checkable. The second case (C2)
can be handled by solving $M\bx=-\bq$. It is the third case (C3) that
needs sophisticated treatments and leads to different numerical methods. For example,  for (C3), \cite{walz:2019,zhys:2015} consider {transforming} $\SOCLCP(\coneK^n,M,\bq)$ into a zero-finding problem. To see this, note that
$$
\eig(MJ_n)=\eig(J_nM)=\eig(M,J_n),
$$
where $\eig(M,J_n)$ denotes the eigenvalues of the matrix pencil $M-\lambda J_n$. In fact,  $M-\lambda J_n=(MJ_n-\lambda I_n)J_n=J_n(J_nM-\lambda I_n)$ implies
$$
\det(M-\lambda J_n)=(-1)^{n-1}\det(MJ_n-\lambda I_n)=(-1)^{n-1}\det(J_nM-\lambda I_n).
$$

\begin{theorem}[{\cite{zhys:2015}}]\label{thm:newlook}
Suppose $s\in \bbR$ is not an eigenvalue of $MJ_n$. Let
\begin{equation}\nonumber\label{eq:xs}
\bx(s)\equiv\begin{bmatrix}
                 x_1(s) \\
                 \bx_2(s) \\
               \end{bmatrix}:=-(M-sJ_n)^{-1}\bq,
\end{equation}
where $J_n$ is given by \eqref{eq:Jn}.
Then $\bx(s)\in \partial(\coneK^n)$ if and only if $x_1(s)>0$ and $h(s)=0$, where
\begin{equation}\label{Sec1_hs}
h(s):=\bx(s)^{\T} J_n\bx(s)=\bq^{\T} (M-sJ_n)^{-\T}J_n(M-sJ_n)^{-1}\bq.
\end{equation}
\end{theorem}

Returning the third case (C3) in Theorem~\ref{thm:3m-e-c}, we have $(M-s_*J_n)\bx=-\bq$ and $\bx\in \partial(\coneK^n)$. Thus
if  $s_*$ is not an eigenvalue of $MJ_n$, then $h(s_*)=0$,
 a zero-finding problem. The problem
can be further simplified whenever  $M$
admits the GUS property \cite{yayu:2013,zhya:2013a,zhya:2014a}. A complete algebraic-geometric characterization
of $\SOCLCP(\coneK^n,M,\bq)$ \eqref{eq:SOCLCP} with such a property has been established in \cite{yayu:2013,zhya:2013a,zhya:2014a}. As one of the key results, it is proved that $MJ_n$ has exactly one positive eigenvalue $\tau$ in $[0,+\infty)$.
 The reader is referred to \cite[Theorem 2.4]{zhys:2015} for more properties.
As a summary, we provide basic facts of the location of the zero root $s_*$  in Corollary \ref{cor:h(s)+-}, Table~\ref{tbl:h(s)} and Figure~\ref{fig:3cases4h(s)}, and suggest \cite{zhys:2015} for a more detailed discussion.

\begin{table}[t]
\renewcommand{\arraystretch}{1.3}
\caption{The sign of $h(s)$ on $(0,\infty)$ in terms of location of $\bq$ \cite{zhys:2015}.}
\begin{center}\label{tbl:h(s)}
{\begin{scriptsize} \tabcolsep0.1in
\begin{tabular}{|c||c|c|} \hline
cases & where is $\bq$?     & solution $x_*$ or $h(s)$ \\ \hline\hline
1 & $\bq\in \coneK^n$ & $ 0=\bx_*\in \SOL(\coneK^n,M,\bq)$\\\hline
2 &  $\bq\in -M\coneK^n$ & $-M^{-1}\bq=\bx_*\in \SOL(\coneK^n,M,\bq)$\\\hline
3 &  $\begin{array}{c}
              \bq\in (-\coneK^n)\backslash(-M\coneK^n),\\
              \bq\not \in \cR(M-\tau J)
            \end{array}$ &
     $h(s)=\left\{\begin{array}{cl}
           - &\mbox{for}\,\,  s\in(0,s_*),\\
           0 &\mbox{for}\,\,  s=s_*,\\
           + &\mbox{for}\,\,  s\in(s_*,\tau),\\
           + &\mbox{for}\,\,  s\in(\tau,\infty).
      \end{array}\right.$\\ \hline
4 &   $\begin{array}{c}
             \bq\in M\coneK^n\backslash \coneK^n,\\
              \bq\not \in \cR(M-\tau J)
            \end{array}$ &
      $h(s)=\left\{\begin{array}{cl}
              + &\mbox{for}\,\,  s\in(0,\tau),\\
              + &\mbox{for}\,\,  s\in(\tau,s_*),\\
              0 &\mbox{for}\,\,  s=s_*,\\
              - &\mbox{for}\,\,  s\in(s_*,\infty).
            \end{array}\right.$ \\ \hline
 5 &  $\begin{array}{c}
             \bq\not\in(-M\coneK^n)\cup M\coneK^n\cup\coneK^n\cup(-\coneK^n),\\
              \bq\not \in \cR(M-\tau J)
            \end{array}$ &
      $h(s)=\left\{\begin{array}{cl}
          - &\mbox{for}\,\,  s\in(0,s_{*;1}),\\
          0 &\mbox{for}\,\,  s=s_{*;1},\\
          + &\mbox{for}\,\,  s\in(s_{*;1},\tau),\\
          + &\mbox{for}\,\,  s\in(\tau,s_{*;2}),\\
          0 &\mbox{for}\,\,  s=s_{*;2}, \\
          - &\mbox{for}\,\,  s\in(s_{*;2},\infty).
           \end{array}\right.$\\ \hline
\end{tabular}
\end{scriptsize}
}\end{center}
\end{table}

\setlength{\unitlength}{0.60mm} 
\begin{figure}[thb]
{\centering
\begin{tabular}{ccc}
\begin{picture}(66,72)(0,0)
\thinlines
\drawvector{0.0}{36.0}{69.0}{1}{0}
\drawvector{10.0}{8.0}{68.0}{0}{1}
\drawdotline{36.0}{76.0}{36.0}{8.0}
\qbezier(6.0,24.0)(31.39,21.79) (34.5,74)
\qbezier(37,74)(36.39,50) (68,40.5)
\drawcenteredtext{40.0}{33.0}{$\tau$}
\drawcenteredtext{25.0}{39.0}{$s_*$}
\drawcenteredtext{8.0}{33.0}{0}
\drawcenteredtext{16.0}{76.0}{$h(s)$}
\drawcenteredtext{36.0}{2.0}{{\rm case 3}}
\end{picture}
&  
\begin{picture}(66,72)
\thinlines
\drawvector{4.0}{36.0}{65.0}{1}{0}
\drawvector{10.0}{8.0}{68.0}{0}{1}
\drawdotline{36.0}{76.0}{36.0}{8.0}
\drawcenteredtext{40.0}{2.0}{{\rm case 4}}
\qbezier(6.0,41.0)(31.39,44) (34.5,74)
\qbezier(37,74)(36.39,49) (68,30.5)
\drawcenteredtext{17.0}{76.0}{$h(s)$}
\drawcenteredtext{40.0}{33.0}{$\tau$}
\drawcenteredtext{8.0}{33.0}{0}
\drawcenteredtext{61.0}{40.0}{$s_*$}
\end{picture}
&
\begin{picture}(66,72)
\thinlines
\drawvector{4.0}{36.0}{65.0}{1}{0}
\drawvector{10.0}{8.0}{68.0}{0}{1}
\drawdotline{36.0}{76.0}{36.0}{8.0}
\drawcenteredtext{40.0}{2.0}{{\rm case 5}}
\qbezier(6.0,24.0)(31.39,25.79) (34.5,74)
\qbezier(37,74)(36.39,49) (68,30.5)
\drawcenteredtext{61.5}{40.0}{$s_{*;2}$}
\drawcenteredtext{22.5}{40.0}{$s_{*;1}$}
\drawcenteredtext{40.0}{33.0}{$\tau$}
\drawcenteredtext{12.0}{33.0}{0}
\drawcenteredtext{17.0}{76.0}{$h(s)$}
\end{picture}
\end{tabular}\par
}
\caption{$h(s)$ corresponding to cases 3, 4, and 5 in Table~\ref{tbl:h(s)} \cite{zhys:2015}.}
\label{fig:3cases4h(s)}
\end{figure}
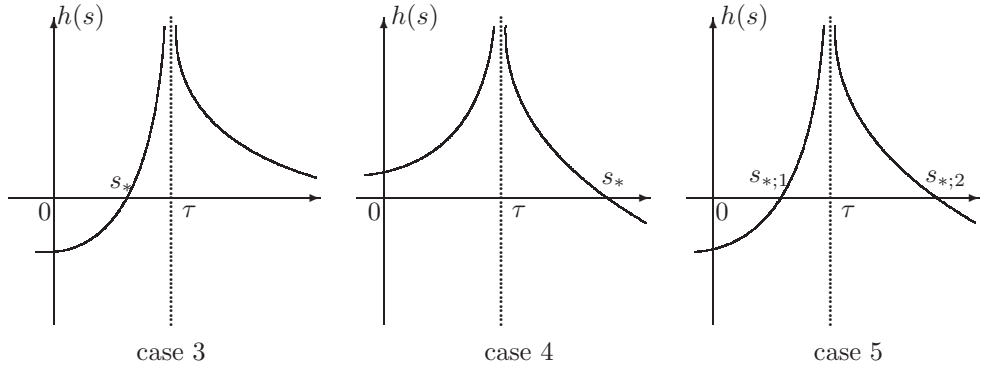

\begin{corollary}\label{cor:h(s)+-}
When $M$ has the GUS property, then the following statements hold:
\begin{enumerate}[{\rm (1)}]
  \item $h(s)$ has a zero in $(0,\tau)$ if and only if $h(0)=\bq^{\T}M^{-\T}J_nM^{-1}\bq<0$;
  \item $h(s)$ has a zero in $(\tau,\infty)$ if and only if $\bq^{\T}J_n\bq<0$.
\end{enumerate}
\end{corollary}

The set of all symmetric and positive definite matrices is a subset of GUS \cite{yayu:2013}. For this special case, we further have more nice properties that can be used for finding the zero root $s_*$.  In particular, for example, according to \cite[Lemma 2.2]{walz:2019} (see also Lemma~\ref{lm:MUJ} later), we know that there is a nonsingular matrix $V\in\bbR^{n\times n}$ such that
\begin{equation}\nonumber\label{eq:eigD-(M,J)}
V^{\T}MV=\Omega\equiv\diag(\omega_1,\omega_2,\ldots,\omega_n),\quad V^{\T}J_nV=J_n,
\end{equation}
where $0<\omega_1=\tau$ and $0<\omega_2\le\cdots\le\omega_n$. In particular,
$\eig(M,J_n)=\{\omega_1,-\omega_i\,\,\mbox{for $2\le i\le n$}\}$.

\section{Rational Krylov Subspace Methods}\label{section:RKSM}
The transformation from $\SOCLCP(\coneK^n,M,\bq)$ to a zero-finding problem $h(s_*)=0$ finishes {\it step 1} within  the numerical framework of \cite{zhys:2015}. We next discuss techniques for {\it step 2} to form proper Krylov subspaces onto which the large-scale $\SOCLCP(\coneK^n,M,\bq)$ can be projected and approximately solved.

\subsection{Projection via Rational Krylov Subspaces}
The motivation of using the rational Krylov subspace is the similarity of the function $h(s)$ \eqref{Sec1_hs} to
transfer functions for certain time-invariant single-input-single-output (SISO)
dynamical systems \cite{anto:2005,sucr:1991}. Both are rational functions that involve the inversions  of
usually parameter dependent matrices, and it is $(M-s J_n)^{-1}$ in our case. In SISO dynamical systems, the transfer function usually has to be evaluated at a wide range of the parameter value, which requires very high computational costs.
Model reduction techniques based on Krylov subspace projections \cite{anto:2005,sucr:1991} are popular and efficient in
mitigating high costs of such evaluations. The basic idea is to generate a suitable Krylov subspace
and then project the transfer function onto the subspace to yield a reduced transfer function, 
 effectively
reduces the original dimension to a few tens or hundreds.
The accuracy of approximation by the reduced transfer function is measured by the number of
leading terms in its Taylor expansion at a point of interest that match
those of the original transfer function in its Taylor expansion at the same point. 
{In model order reduction area, rational Krylov subspace methods are demonstrated efficient in computing system's {$\cal{L}_\infty$} norm \cite{MR3735291} or  pseudospectral abscissa  \cite{MR3180856}, and  solving algebraic {R}iccati equations \cite{MR3232439}.
For more details, the reader is referred to, e.g.,
\cite{anto:2005,basu:2005a,fefr:1995,gagv:1994,liba:2005,liye:2014,scvr:2008,sucr:1991}
and references therein.

Inspired by the idea in the model reduction, \cite{zhys:2015}
proposes a Krylov subspace method to find the zero $s_*$
of $h(s)$. It is an iterative method that can roughly be explained as follows.
Suppose $0<s_0\in \mathbb{R}$ is a given approximation to $s_*$ and let
\begin{equation}\nonumber\label{eq:shiftM}
A_{s_0}=(M-s_0J_n)^{-1}J_n, \quad
\bb_{s_0}=-(M-s_0J_n)^{-1}\bq.
\end{equation}
We have $\bx(s)=-(M-sJ_n)^{-1}\bq= [I_n-(s-s_0) A_{s_0}]^{-1} \bb_{s_0}$, and
\begin{equation}\nonumber\label{eq:h(s)-shifted}
h(s)=\bx(s)^{\T} J_n  \bx(s)
    =\bb_{s_0}^{\T} [I_n-(s-s_0) A_{s_0}]^{-\T} J_n [I_n-(s-s_0) A_{s_0}]^{-1} \bb_{s_0}.
\end{equation}
Suppose $s_0$ is neither a pole nor a zero of $h(s)$ (in other words, $h(s_0)\ne 0$ and is finite).
This implies $h(s_0)=\bb_{s_0}^{\T} J_n\bb_{s_0}\neq 0$. Let $Y_{\ell}\in\bbR^{n\times \ell}$ be
an orthonormal basis matrix of the Krylov
subspace $\cK_{\ell}(A_{s_0},\bb_{s_0})$, i.e.,
$Y_{\ell}^{\T}Y_{\ell}=I_{\ell}$ and $\cR(Y_{\ell})=\cK_{\ell}(A_{s_0},\bb_{s_0})$. Then the reduced $h(s)$ is given by
\begin{equation}\nonumber\label{eq:h(s)-reduced:zhys2005}
h_{\ell}(s)=\|\bb_{s_0}\|_2^2\be_1^{\T} [I_{\ell} -(s-s_0) H_{\ell}]^{-\T}
             Y_{\ell}^{\T} J_n Y_{\ell}
             [I_{\ell} -(s-s_0) H_{\ell}]^{-1}\be_1,
\end{equation}
where $H_{\ell}=Y_{\ell}^{\T}A_{s_0}Y_{\ell}$. Both $Y_{\ell}$ and $H_{\ell}$ can be efficiently computed by
the Arnoldi Process \cite[Algorithm 1]{zhys:2015}.
It is proved that \cite{zhys:2015}
\begin{equation}\label{eq:h(s)-mm}
h(s)=h_{\ell}(s)+{O}(|s-s_0|^\ell).
\end{equation}

Indeed, an immediate implication of \eqref{eq:h(s)-mm} is that the first $\ell$ leading terms
of the Taylor expansions of $h(s)$ and $h_{\ell}(s)$ at $s_0$ match, or equivalently,
the coefficients of $(s-s_0)^i$ for $i=0,1,\ldots,\ell$, called { moments} in the expansions,
are the same.
A particular zero of $h_{\ell}(s)$ is then computed as the next approximation.
Since $h_{\ell}(s)$ in general has many zeros, with the help of Table~\ref{tbl:h(s)} and Figure~\ref{fig:3cases4h(s)},
the method of \cite{zhys:2015} picks a particular positive  one; conceivably, if $h_{\ell}(s)$ approximates $h(s)$ sufficiently well in
the region of interest, a positive zero root exists and solving $h_{\ell}(s)=0$ can be done by  equivalently transforming it into an $(\ell-1)\times(\ell-1)$ quadratic eigenvalue problem; the process repeats whenever the updated approximation is not within the given accuracy.

Our improvement of the rational Krylov subspace method over \cite{zhys:2015}  begins with the following approximation theorem when $M \succ 0$.
\begin{theorem}\label{thm:h(s)-shifted}
If $M$ is symmetric, then
\begin{equation}\label{eq:h(s)-mm'}
h(s)=h_{\ell}(s)+{O}(|s-s_0|^{2\ell-1}).
\end{equation}
\end{theorem}

\begin{proof}
Let $f(s)=q^{\T}(M-sJ)^{-1}q$ and $f_{\ell}(s)=\|\bb_{s_0}\|_2^2\be_1^{\T} [I_{\ell} -(s-s_0) H_{\ell}]^{-\T}\be_1$.
It can be verified that
\begin{equation}\label{eq:h(s)-shifted:pf-1}
f'(s)=h'(s), \quad f_{\ell}'(s)=h_{\ell}'(s).
\end{equation}
By \cite[Theorem 3.3]{liba:2005}, we find $f(s)=f_{\ell}(s)+{O}(|s-s_0|^{2\ell})$.
The equation \eqref{eq:h(s)-mm'} is a consequence of \eqref{eq:h(s)-shifted:pf-1}.
\end{proof}

We remark that the approximation accuracy as measured by moment matching  given in \eqref{eq:h(s)-mm} is the best possible
when $M$ is in general non-symmetric. For a symmetric $M$,  we see from Theorem \ref{thm:h(s)-shifted} that the number of matched moments doubles.

Instead of \cite{zhys:2015}'s subspace
$$
\cK_{\ell}(A_{s_0},\bb_{s_0})\equiv\cK_{\ell}((M-s_0J_n)^{-1}J_n,(M-s_0J_n)^{-1}\bq),
$$
we propose to  use a new
rational Krylov subspaces to reduce $h(s)$, namely,  the sum of several Krylov subspaces in the form of
$\cK_{\ell}(A_{s_i},\bb_{s_i})$,  each of which is expanded at a different point $s_i$. Specifically, let $U$ be the orthonormal
matrix such that
\begin{equation}\label{eq:sum_of_rational_Krylov}
\cR(U)=\sum_{i=1}^j\cK_{\ell}((M-s_iJ_n)^{-1}J_n,(M-s_iJ_n)^{-1}\bq),
\end{equation}
and we then reduce $h(s)$ to
\begin{subequations}\label{eq:h(s)-RKS}
\begin{equation}\label{eq:h(s)-RKS:1}
\hat h(s)=\hat\bq^{\T}( \what M-s \what J)^{-\T} \what J( \what M-s \what J)^{-1}\hat\bq,
\end{equation}
 where
\begin{equation}\label{eq:h(s)-RKS:2}
\what M=U^{\T}MU,\quad \what J=U^{\T}J_n U, \quad \hat\bq=U^{\T}\bq.
\end{equation}
\end{subequations}
By using various  $s_i$ around $s_*$, it is hoped that
the reduced function $\hat h(s)$ is able to qualitatively match better $h(s)$, and therefore, the  overall convergence behavior can be improved.
%


\subsection{Solve $\hat h(s)=0$}\label{ssec:th(s)=0}
Next, we consider {\it step 3} to solve the zero of the  reduced $\hat h(s)=0$. We point out that the development in this subsection for $\hat h(s)=0$ of form \eqref{eq:h(s)-RKS} works for
any $U\in\bbR^{n\times m}$ with full column rank, i.e., $\rank(U)=m$, not restricted in the basis  from a rational Krylov subspace.


To utilize $ \hat h(s) \approx h(s)$, we need the assumption}
\begin{equation}\label{assume:UtJU}
\framebox{
\parbox{12cm}{
$U^{\T}J_nU$ has one positive eigenvalue and the rest of its eigenvalues are negative.
}
}
\end{equation}
When $U$ has orthonormal columns,  this should hold
as the number of columns of $U$ increases. 
In general,   $U^{\T}J_nU$ has at most one positive eigenvalue.

\begin{lemma}\label{lm:UJU}
Suppose that $U\in\bbR^{n\times m}$ has orthonormal columns, where $1\le m<n$. Then $U^TJ_nU$ has at most one positive eigenvalue.
\end{lemma}

\begin{proof}
Let $U_{\bot}\in\bbR^{n\times (n-m)}$ such that $[U,U_{\bot}]$ is orthogonal. The eigenvalues of
$$
[U,U_{\bot}]^{\T}J_n[U,U_{\bot}]=\begin{bmatrix}
                                U^{\T}J_nU & U^{\T}J_nU_{\bot} \\
                                U_{\bot}^{\T}J_nU & U_{\bot}^{\T}J_nU_{\bot}
                              \end{bmatrix}
$$
are the same as $J_n$. Denote by $\mu_1\ge\mu_2\ge\cdots\ge\mu_m$ the eigenvalues of $U^{\T}J_nU$. It suffices
to show all $\mu_i<0$ for $2\le i\le m$. By the Cauchy's interlacing inequalities \cite{hojo:2013}, we have
$$
1\ge\mu_1\ge -1,\quad -1\ge\mu_i\ge -1\,\,\mbox{for $2\le i\le m$}
$$
implying $\mu_i=-1<0$ for $2\le i\le m$, as expected.
\end{proof}

%
%
%
%

Now, we describe two methods for solving $\hat h(s)=0$: the first one closely follows
the idea in \cite{zhys:2015} by turning it into a quadratic eigenvalue problem, and it works for both symmetric and nonsymmetric $M$; the second one
is the method of \cite{walz:2019} and it works for $M\succ 0$ only.

For the first method, we note that the direct transformation to a quadratic eigenvalue problem \cite{zhys:2015}  should be modified to work on \eqref{eq:h(s)-RKS}.  In particular, suppose  $\hat\bq\in\bbR^m$ and  $Q\in\bbR^{m\times m}$ is an orthogonal matrix (for example, the Householder matrix \cite{demm:1997}) such that
$
Q^T\hat\bq=\|\hat\bq\|_2\be_1.
$
We then can write $\hat h(s)$ as follows:
$$
\hat h(s)=\|\hat\bq\|_2^2\,\be_1^{\T}\Big[s^2Q^{\T}\what JQ-sQ^{\T}(\what M^{\T}+\what M)Q+Q^{\T}\what M\what J^{-1}\what M^{\T}Q\Big]^{-1}\be_1.
$$
With this preprocess, the method of \cite[subsection 6.2]{zhys:2015} then is applicable.

For the second method related with \cite{walz:2019}, we need the next lemma which is a straightforward extension of
\cite[Lemma 2.2]{walz:2019} and can be proved in a similar way.

\begin{lemma}\label{lm:MUJ}
Suppose $M \succ 0$ and assume \eqref{assume:UtJU} holds. Then there is a nonsingular matrix $\what V\in\bbR^{m\times m}$ such that
\begin{equation}\label{eq:eigD-(M,UJU)}
\what V^{\T}U^{\T}MU\what V=\what\Omega\equiv\diag(\hat\omega_1,\hat\omega_2,\ldots,\hat\omega_m),\quad \what V^{\T}U^{\T}J_nU\what V=J_m,
\end{equation}
where $0<\hat\omega_1$ and $0<\hat\omega_2\le\cdots\le\hat\omega_m$.
\end{lemma}

\begin{proof}
Write $\what M=U^{\T}MU$ and $\what J_m=U^{\T}J_nU$.
Since $\what M\succ 0$, it has a Cholesky decomposition $\what M=R^{\T}R$.
Now notice that $R^{-\T}\what J_mR^{-1}\in\bbR^{m\times m}$ is symmetric and let its eigenvalues
be $\{\mu_i\}_{i=1}^m$. Because $R^{-\T}\what J_mR^{-1}\in\bbR^{m\times m}$  has
the same inertia as $J_m$,  these eigenvalues can be ordered in such a way that
\begin{equation}\label{eq:eigD-(M,UJU):pf-0}
\mu_m\le\cdots\le\mu_2<0<\mu_1.
\end{equation}
$R^{-\T}\what J_mR^{-1}$ has an eigendecomposition
\begin{equation}\nonumber\label{eq:eigD-(M,UJU):pf-1}
R^{-\T}\what J_mR^{-1}=U D U^{\T},\quad D=\diag(\mu_1,\mu_2,\ldots,\mu_m),
\end{equation}
where $U$ is an orthogonal matrix. Set $\hat\omega_i=1/|\mu_i|$ for $1\le i\le n$. We have
$$
R^{-\T}\what J_mR^{-1}=U D U^{\T}=UJ_m\what\Omega^{-1}U^{\T}=U\what\Omega^{-1/2}J_m\what\Omega^{-1/2}U^{\T}.
$$
Finally set $\what V=R^{-1}U\what\Omega^{1/2}$ to conclude the proof.
\end{proof}

Using the decompositions in \eqref{eq:eigD-(M,UJU)}, we have
\begin{align}\nonumber
\hat h(s)&={\wtd\bq}^{\T}(\what\Omega-s J_m)^{-1}J_m(\what\Omega-s J_m)^{-1}\wtd\bq \nonumber\\
    &=\frac {\xi_1^2}{(s-\hat\omega_1)^2}-\sum_{i=2}^m\frac {\xi_i^2}{(s+\hat\omega_i)^2},
      \label{eq:th(s)-explicit}
\end{align}
where $\xi_i$ for $1\le i\le m$ are the entries of $\what V^{\T}U^{\T}\bq$, i.e.,
\begin{equation}\nonumber\label{eq:VTUTq}
\wtd\bq=\what V^{\T}\hat\bq=\what V^{\T}U^{\T}\bq=[\xi_1,\ldots,\xi_m]^{\T}.
\end{equation}
This $\hat h(s)$ takes the same form as the one of \cite[(3.1)]{walz:2019}, and has up to
two positive zeros as indicated by Figure~\ref{fig:3cases4h(s)}.
{Note that $\hat\omega_1$ plays the similar pole for $\hat h(s)$ as $\tau$ does for $ h(s)$.} We can use the efficient zero-finder
there to find the zeros.


\begin{remark}\label{rk:key-condition}
In the process of forming $\hat h(s)$ in  \eqref{eq:th(s)-explicit}, as a by-product, we have also proved that the condition \eqref{assume:UtJU} holds.
In fact, \eqref{assume:UtJU} is equivalent to \eqref{eq:eigD-(M,UJU):pf-0},
which is used in our implementation.
\end{remark}

There is a subtle but important comment to make. 
Suppose $h(s)=0$ of $\SOCLCP(\coneK^n,M,\bq)$ has two positive zeros $s_{*;1}<s_{*;2}$ one of which corresponds to the solution (cf. Figure~\ref{fig:3cases4h(s)}). 
Associated with $\hat h(s)=0$ is a reduced $\SOCLCP(\coneK^m,\what\Omega,\wtd\bq)$ with $\what\Omega\succ 0$, and  $\hat h(s)$ may also have two positive zeros, say $\hat s_{*;1}<\hat s_{*;2}$, and suppose $\hat s_{*;1}$ corresponds to  the solution of $\SOCLCP(\coneK^m,\what\Omega,\wtd\bq)$. However, we cannot ensure $s_{*;1}$ is the right one for the solution of $\SOCLCP(\coneK^n,M,\bq)$ (cf. Figure~\ref{fig:3cases4h(s)}).
This is observed in our numerical results. For that reason, we cannot simply exclude one of two zeros in $\hat h(s)=0$ based on $\SOCLCP(\coneK^m,\what\Omega,\wtd\bq)$, and thus have to compute both zeros. A more detailed procedure for determining  the right zero in $\hat h(s)=0$ will be given in Figure \ref{fig:philosophy_of_case345}.
\subsection{Initial approximation subspace}\label{ssec:initial}
Though Figure~\ref{fig:3cases4h(s)} clearly gives the location of $s_*$ relative to $\tau$, the only positive eigenvalue of
$M J_n$, it is still costly for computing $\tau$ when $n$ is large. 
Fortunately, we do not need to compute $\tau$.
The algorithm is expected to start from an initial subspace, which can provide enough information, and are not 
expensively constructed.

Since the classical Lanczos method is able to approximate well  the extreme eigenvalues of a square matrix \cite{zhli:2011}, and $\tau$ is the only positive eigenvalue of $MJ_n$ (also of $J_nM$),  we propose to
initially build an extended Krylov subspace \cite{knsi:2011,simo:2007}
of $J_nM$ on $J_n\bq$ as
\begin{equation}\nonumber\label{eq:extK-initial}
\cK_{\ell_0,k_0}^{\ext}(J_nM,J_n\bq)=\cK_{\ell_0}(J_nM,J_n\bq)+\cK_{k_0}((J_nM)^{-1},(J_nM)^{-1}J_n\bq).
\end{equation}
Let $U$ be an orthonormal basis of $\cK_{\ell_0,k_0}^{\ext}(J_nM,J_n\bq)$, i.e.,
$U^{\T}U=I$ and $\cR(U)=\cK_{\ell_0,k_0}^{\ext}(J_nM,J_n\bq)$. We then form $U^{\T}MU$, $U^{\T}\bq$, and $U^{\T}J_nU$ to give
$\hat h(s)$ of \eqref{eq:h(s)-RKS}.
{Note $\cK_{k_0}((J_nM)^{-1},(J_nM)^{-1}J_n\bq)$
provides shift $s_0=0$ in \eqref{eq:h(s)-mm'}.
Thus, we get $h(0)=\hat h(0).$}

There are two scenarios for the solution $\hat h(s)=0$. Case (1): when \eqref{assume:UtJU} is true, then the method described in subsection~\ref{ssec:th(s)=0} is able to find the particular zero of $\hat h(s)=0$, 
which provides the next shift $s_j$. Case (2):  when \eqref{assume:UtJU} fails, it suggests that the extended Krylov subspace is not big enough; as an economic treatment for the latter, by a fact  $\tau\le\|M\|_1$, we choose to set the approximations $s_j=\|M\|_1/10^{j}$ { and add $\cK_{\ell}((M-s_jJ_n)^{-1}J_n,(M-s_jJ_n)^{-1}\bq)$ to the obtained space}.

\subsection{The main algorithm}


With the initial extended Krylov subspace $\cK_{\ell_0,k_0}^{\ext}(J_nM,J_n\bq)$, we now describe our  rational Krylov subspace methods (RKSM) for $\SOCLCP$ \eqref{eq:SOCLCP}.

The main procedure of RKSM is to expand the rational Krylov subspace. Suppose that we have already generated $j$ approximations  $s_1,\ldots,s_j$, together with an orthonormal basis matrix $U$ of the rational Krylov subspace
of \eqref{eq:sum_of_rational_Krylov}. 
We solve the reduced problem $\hat h(s)=0$ for the next $s_{j+1},$ and add $$
\cK_{\ell}((M-s_{j+1}J_n)^{-1}J_n,(M-s_{j+1}J_n)^{-1}\bq)
$$ to obtain a new $U$.

 As we remarked at the end of  Section \ref{ssec:th(s)=0}, computing the right zero of each reduced $\hat h(s)=0$ should be carefully treated. Let $\hat \omega_1$ be given by \eqref{eq:eigD-(M,UJU)}. Based on the facts revealed in Figure \ref{fig:3cases4h(s)}, we  use the procedure in Figure~\ref{fig:philosophy_of_case345}. It is similar to the one in  \cite[Algorithm 5.1]{walz:2019}.

\begin{figure}[hbt]
%
%
%
    \begin{algorithmic}[1]
\IF{$\hat h(0)<0$}
\STATE  compute the zero root $s_{*}$ of $\hat h(s)=0$  in $ (0,\hat\omega_1)$; \qquad \%  case 3 or left of case 5
\STATE $\bx=-(M-s_{*}J)^{-1}\bq;$
\STATE {\bf if} $\bx \in \partial(\coneK^n)$ {\bf then}  {\bf return};  {\bf end if}
\ENDIF
\STATE   compute the zero root $s_{*}$ of $\hat h(s)=0$ in $(\hat\omega_1,\infty)$; \qquad \quad  \%  case 4 or right of case 5
\STATE $\bx=-(M-s_{*}J)^{-1}\bq;$
\STATE {\bf if} $\bx \in \partial(\coneK^n)$ {\bf then}  {\bf return};  {\bf end if}
    \end{algorithmic}
 \caption{{ The strategy of finding zero(s) of $\hat h(s)=0.$ \label{fig:philosophy_of_case345}}}
\end{figure}
 
Finally, we outline RKSM in Algorithm~\ref{alg:RKS4SOCLCP}. {Note the eigenvalue $\tau$ is never computed in the algorithm.}

\begin{algorithm}
\caption{RKSM for $\SOCLCP(\coneK^n,M,\bq)$} \label{alg:RKS4SOCLCP}
\begin{algorithmic}[1]
\REQUIRE $\bq\in\bbR^n$, $M\in \bbR^{n\times n}$ (symmetric positive definite);
\ENSURE  $\bx\in \SOL(\coneK^n,M,\bq)$.

\STATE {\bf if}  $\bq\in \coneK^n$ {\bf then} $\bx= 0$, {\bf return};
\STATE {\bf if} $-M^{-1}\bq\in \coneK^n$ {\bf then} $\bx=-M^{-1}\bq$, {\bf return};

\STATE form an orthonormal basis matrix of an initial subspace, $U=\orth(\cK_{\ell_0,k_0}^{\ext}(J_nM,J_n\bq))$;

\STATE  $\what M=U^{\T}MU$, $\what J=U^{\T}J_nU$, $\hat\bq=U^{\T}\bq$;


\IF{$h(0)<0$}
\FOR{$j=1,2, \ldots ,j_{\max}$}

\IF{$\what J$ has no positive eigenvalue}
     \STATE $s_j:=\frac{\|M\|_1+\texttt{fix}(\frac{k}{16})}{10^{\texttt{mod}(j,16)}}$;
\ELSE
     \STATE compute the (smaller) positive zero of $\hat h(s)$, called it $s_j$;
\ENDIF

\STATE   $\bx=-(M-s_jJ_n)^{-1}\bq;$
\STATE  $h=\bx^{\T}J_n\bx;$ {\bf if} $|h|<\epsilon_1\,\|\bx\|_2^2$ {\bf then} {\bf break};  {\bf end if}
\STATE  {\bf if} $\chi_{\rel}(\bx)<\epsilon_2$ {\bf then} {\bf return};  {\bf end if}
 \STATE $V=\orth(\cK_{\ell}((M-s_jJ_n)^{-1}J_n,(M-s_jJ_n)^{-1}\bq))$;
          \STATE $U=\orth([U,V])$;
     \STATE $\what M=U^{\T}MU$, $\what J=U^{\T}J_nU$, $\hat\bq:=U^{\T}\bq$; 
\ENDFOR
\STATE {\bf if} $\bx\in \partial(\coneK^n)$ {\bf then}  {\bf return};  {\bf end if}
\ENDIF
%
\FOR{$j=1,2, \ldots ,j_{\max}$}

\IF{$\what J$ has no positive eigenvalue}
     \STATE $s_j=(1.1)^{j-1} \times\|M\|_1$;
\ELSE
        \STATE compute the (larger) positive zero of $\hat h(s)$,  called it $s_j$;
\ENDIF

\STATE $\bx=-(M-s_jJ_n)^{-1}\bq;$
\STATE  $h=\bx^{\T}J_n\bx;$ {\bf if} $|h|<\epsilon_1\,\|\bx\|_2^2$ {\bf then} {\bf break};  {\bf end if}
\STATE  {\bf if} $\chi_{\rel}(\bx)<\epsilon_2$ {\bf then} {\bf return};  {\bf end if}
 \STATE $V=\orth(\cK_{\ell}((M-s_jJ_n)^{-1}J_n,(M-s_jJ_n)^{-1}\bq))$;
      \STATE $U=\orth([U,V])$;
     \STATE $\what M=U^{\T}MU$, $\what J=U^{\T}J_nU$, $\hat\bq=U^{\T}\bq$; 
\ENDFOR
%

\IF{$\bx \in \partial(\coneK^n)$}
\STATE  {\bf return};
\ELSE
\STATE Call \cite[Algorithm 2]{zhya:2013a} for the special case: $s_*$ is the positive eigenvalue of $MJ_n$.
\ENDIF
\end{algorithmic}
\end{algorithm}

\begin{remark}\label{rk:alg:RKS4SOCLCP}
We now provide implementation details for Algorithm~\ref{alg:RKS4SOCLCP} and make comments.
\begin{enumerate}
  \item {\bf Lines 1 and 2.} Check Cases 1 and 2 in Table~\ref{tbl:h(s)}. The Cholesky  factor of $M$ at Line 2 can be used for the following computations.

  \item {\bf Line 3.} Build orthonormal bases $U_1$ and $U_2$ for
        $\cK_{\ell_0}(J_nM,J_n\bq)$ and $\cK_{k_0}((J_nM)^{-1},(J_nM)^{-1}J_n\bq)$,
       respectively, by the Arnoldi process.
The Cholesky factor of $M$ can be used for  $U_2$. The orthonormalization for $U_2$ against $U_1$ is as follows
 $$
       U_2\leftarrow U_2-U_1(U_1^{\T}U_2), \quad U=[U_1,\orth(U_2)].
       $$

An efficient procedure  \cite{simo:2007} can be employed  for $k_0=\ell_0$ to form the orthonormal basis matrix.
     \item {\bf Line 5.}
         $h(0)=\bq^{\T}M^{-\T}J_nM^{-1}\bq$ is computed at marginal cost with $M^{-1}\bq$ from Line 2.
       Here, we  check whether $h(0)<0$  or not to conclude
       if $h(s)$ has a positive zero in $(0,\tau)$ (Corollary~\ref{cor:h(s)+-}).
 Corresponding to the two cases, the for-loop from Line 6 to 18  computes the zero point in $(0,\tau)$, while the for-loop from Line 21 to 33  computes the zero point in $(\tau,\infty)$.

   \item {\bf Lines 6 and 21.} We set $J_{\max}=40$. If $h(0)<0$ at Line 5 and the algorithm doest not break out of the for-loop from Line 6 to 18, it indicates a failure in finding a  zero point within the maximal number of iterations $J_{\max}=40$. The similar statement applies to Line 21.

  \item {\bf Lines 7 and 22.} Two scenarios of the reduced system $\hat h(s)=0$ (see  Remark~\ref{rk:key-condition}).

  \item{\bf Line 8.} Indicate that the reduced problem cannot generate positive the shifts, and then use shift $s_j=\|M\|_1/10^{j}$ instead. We choose $16$ as the maximal number of iterations because the shift in this case is closed to Matlab \texttt{eps}, for which the associated subspace becomes ineffectively. New shifts will be used after $16$ iterations.
%

  \item {\bf Lines 10 and 25.}    Compute the positive zero by the method described in  subsection~\ref{ssec:th(s)=0}. The smaller (larger) one is selected if there exist two zero points.

  \item {\bf Lines 12 and 27.} The LDL decomposition of $M-s_jJ_n$ are  kept for reuse at Lines 15 and 30 in building an orthonormal basis matrix $V$ by the Arnoldi process.

  \item {\bf Lines 13 and 28.} Stopping criteria for the zero-finding problem in the for-loop ($\epsilon_1=10^{-7}$ in our testing). This rule is based on the rounding error  fact
  $$
        |\fl(\bx^{\T} J_n\bx)-\bx^{\T}J_n\bx|\le 2n\macheps \|\bx\|_2^2,
        $$
        where $\fl(\bx^{\T}J_n\bx)$ denotes the computed $\bx^{\T}J_n\bx$, given $\bx$, and $\macheps$ is the unit machine roundoff.  

\item {\bf Lines 14, and 29.} Check the relative error $\chi_{\rel}(\bx)<\epsilon_2$ (default $10^{-8}$ )  to stop the iteration.
The total relative error $\chi_{\rel}(\bx)$ is defined by \cite{zhys:2015}
\begin{flalign}
\chi_{\rel}(\bx)&=\chi_{\rel_1}+\chi_{\rel_2}+\chi_{\rel_3},\label{eqn:relative_error}\\
\chi_{\rel_1}&=\frac {\max\{\|\bx_{(2:n)}\|_2-x_1,0\}}{\|\bx\|_2},\nonumber \\
\chi_{\rel_2}&=\frac {\max\{\|\bg_{(2:n)}\|_2-g_1,0\}}{\|M\|_1\|\bx\|_2+\|\bq\|_2},\nonumber \\
\chi_{\rel_3}&=\frac {|\bx^{\T} \bg|}{\|\bx\|_2(\|M\|_1\|\bx\|_2+\|\bq\|_2)},\nonumber
\end{flalign}
where $\bg=M\bx+\bq$. The rationality of these relative errors is explained in \cite{zhys:2015}.
The cost of computing $\chi_{\rel}(\bx)$ is affordable.

  \item {\bf Lines 16 and 31.} Compute orthonormal basis matrix for the combined subspace by       $$
       V\leftarrow V-U(U^{\T}V), \quad U=[U,\orth(V)].
       $$

  %

\item {\bf Lines 19  and 34.} Check the solution for (C3) of Theorem \ref{thm:3m-e-c} on $\partial(\coneK^n)$.
   In particular, check if $x_1>0$ and $|x_1-\|\bx_{(2:n)}\|_2|<\epsilon_3\,\|\bx\|_2$  with defaulting setting $\epsilon_3=10^{-6}.$

  \item{\bf Line 23.} This can be executed only if  $h(0)>0$. At this moment (case 4 in Figure~\ref{fig:3cases4h(s)}), in order to span a more effective subspace, we introduce  a shift $s_j:=(1.1)^{j-1} \times\|M\|_1$. It may happen that $h(s_{\check{j}})=\hat h(s_{\check{j}})<0$ for some $\check{j}$. In this case, because $\hat h(s)>0$ for $s$ near $\omega_1$,  there is a zero point in $(\omega_1,s_{\check{j}})\subset (0,s_{\check{j}})$. Thus, we can compute the positive  point of $\hat h(s)=0$ in that interval. By Remark~\ref{rk:key-condition}, we  know that the condition at {Line 22}  will be true after we introduce some compulsory shifts.

%
    \item{\bf Line 37.} The treatment \cite[Algorithm 2]{zhya:2013a} for the special case (i.e. the assumption of Theorem \ref{thm:newlook} does not hold). 

\end{enumerate}
\end{remark}

Finally, we remark that RKSM differs from   LCPvA \cite{zhys:2015} in the following { three } aspects:
\begin{itemize}
\item [(1)] We use an initial subspace (an extended Krylov subsapce) to obtain the approximates of the zero points, which can fasten the convergence;
\item[(2)] We accumulate all the computed subspaces bases, while LCPvA discards the previous subspaces bases. This could make RKSM more efficient and robust;
\item[(3)]  {Our strategy for solving $h(s)=0$ follows the one used in \cite[Algorithm 5.1]{walz:2019}: we first compute the zero point in $(0,\tau)$ and check whether it is the solution; if it is not, we do the similar process in $(\tau,\infty)$. LCPvA computes approximations  of all  zero roots, and choose one  satisfying the stopping criterion $\chi_{\rel}(\bx) <10^{-7}$ (cf. (\ref{eqn:relative_error})).}

\end{itemize}
\section{Numerical experiments}\label{section:numerical_experiments}
This section is devoted to the evaluation of our proposed RKSM. We carry out our numerical testings upon the Matlab2016a  platform on a notebook (64 bits) with an Intel CPU i7-5500U and  8GB memory. The tests use a fixed \texttt{q=ones(n,1)} and various types of $M$.

\subsection{Typical behavior of RKSM}
We first use the following matrices to see the performance of RKSM.

\vskip 3mm

\noindent\textbf{Example 1.}  Set
\begin{equation}\label{eq:sprandsym}
R={\rm \texttt{sprandsym(n,density,rc},\texttt{kind)}},
\end{equation}
with $\texttt{n}=3000, \texttt{desity}=0.005,\texttt{rc}=0.01,\texttt{kind}=2$ and  $M=R^{\T}R\succ 0.$

Choose  $\cK_{3,3}^{\ext}(J_nM,J_n\bq)$ as initial subspace at line 3 in Algorithm \ref{alg:RKS4SOCLCP}. Note $$\dim(\cK_{3,3}^{\ext}(J_nM,J_n\bq))=6.$$ 
At lines $15$ and $30$, we let $\ell=1$, and set the stopping criterion at lines $13$ and $28$ as
$\epsilon_1=10^{-8}$. To obtain a relatively high accurate solution, the stopping criterion at lines $14$ and $29$ is  $\epsilon_2=10^{-12}$.

With a particular case of { Example 1}, we observed that  the testing problem falls into case 5 in Figure \ref{fig:3cases4h(s)}. Also, the positive eigenvalue of $MJ$ is $\tau= 0.107571187229409$. 
For the shifts used during the iteration, we list them in Table \ref{table:all_shifts}.
Table \ref{table:convergence_order} gives the quantities of $|s_j-s_{*;1}|$ and $\frac{|s_j-s_{*;1}|}{|s_{j-1}-s_{*;1}|^2}$  which reflects the order of the quadratic  convergence.
Note that both $s_1$ and $s_2$ are larger than $\tau$,
but are not close enough to $s_{*;1}$.

\begin{table}[htbp]
\centering
\caption { All of shifts of Example 1 \label{table:all_shifts}.}
\begin{tabular}{|l|l|l|l|l|l|}
\hline
initial subspace &  \multicolumn{5}{c|}{$\cK_{3,3}^{\ext}(J_nM,J_n\bq)$} \\ \hline
 shifts  $s_j$ at Line 10  & 1.3520 & 0.1997& 0.1043 & 0.1034 & 0.1034  \\ \hline
 shifts $s_j$ at Line 25  & \multicolumn{5}{c|}{0.111859222882879}         \\ \hline
\end{tabular}
\end{table}


\begin{table}[htbp]
 \centering
 \caption { Shifts $s_j$ at line 10 converge to $s_{*;1}$ in Example 1\label{table:convergence_order}.}
  \begin{tabular}{ |c|c|c|c| }
\hline
iteration  & shifts $s_j$ &$|s_j-s_{*;1}|$&$ \frac{|s_j-s_{*;1}|}{|s_{j-1}-s_{*;1}|^2} $ \\
\hline
   $1$  & $1.352032867762176 $&  $1.2486e+00 $   & $  -   $   \\
   $2$  & $ 0.199685173922384$  & $9.6238e-02$   &$  6.1732e-02$\\
   $3$  &$ 0.104251405957903 $ &$ 8.0464e-04$  & $8.6877e-02 $\\
   $4 $  &$0.103446765760274$& $   3.8392e-09 $  &$5.9298e-03 $\\
   $5 $ & $ 0.103446769593154$ &  $0$  & $0 $\\
\hline

  \end{tabular}
 \\
 \small{ We use $s_{*;1}= s_5$.}
\end{table}

%
%

\begin{figure}[htb]
  \caption{Behaviors of RKSM for Example 1.}
  \label{fig:kind2_1000_rc0.01} 
%
      \subfigure[$h(s)$]{
    \includegraphics[width=2.9in]{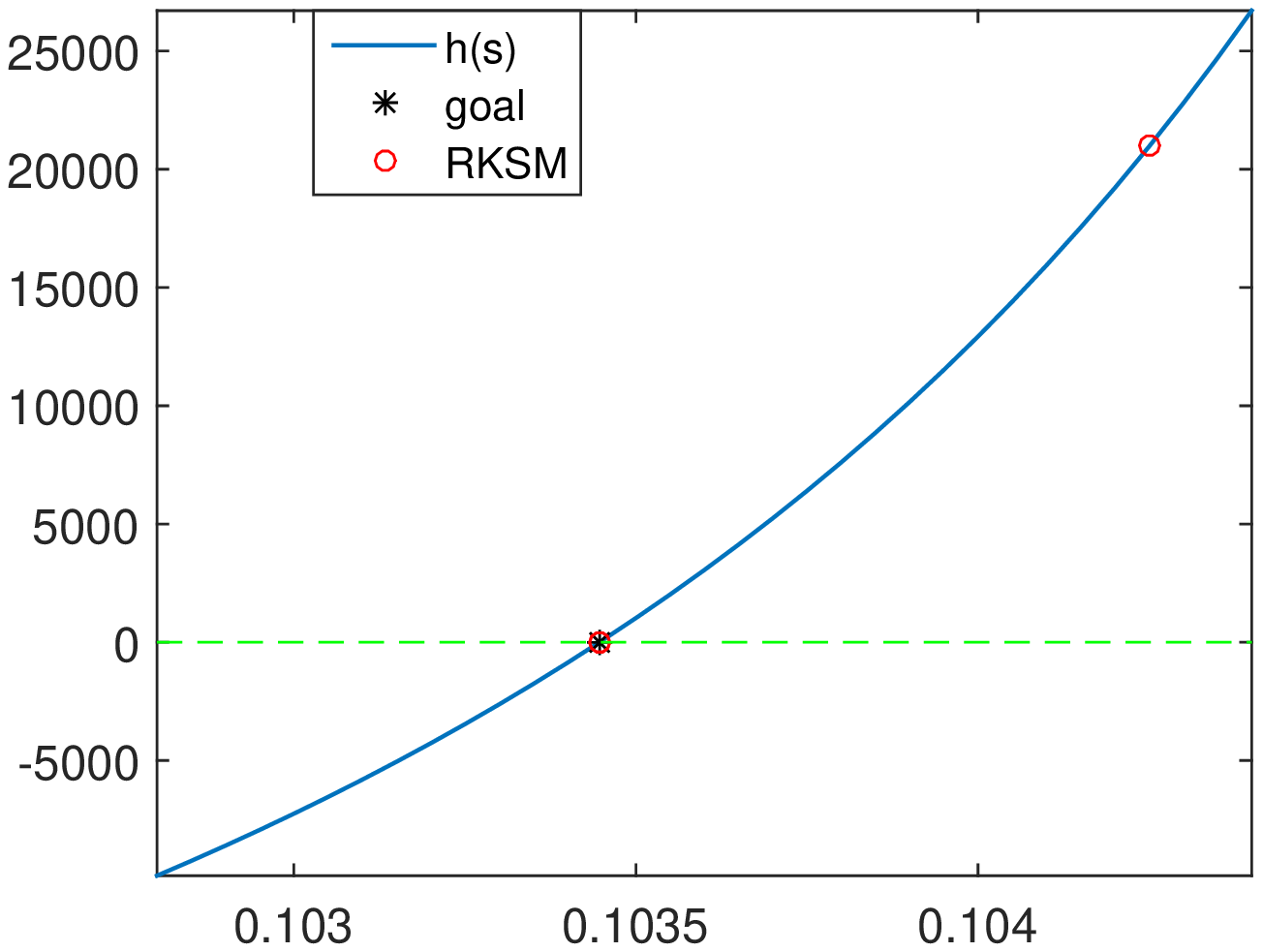}}
   \subfigure[$\frac{|h(s_j)|}{\|x\|^2_2}$]{
    \includegraphics[width=2.9in]{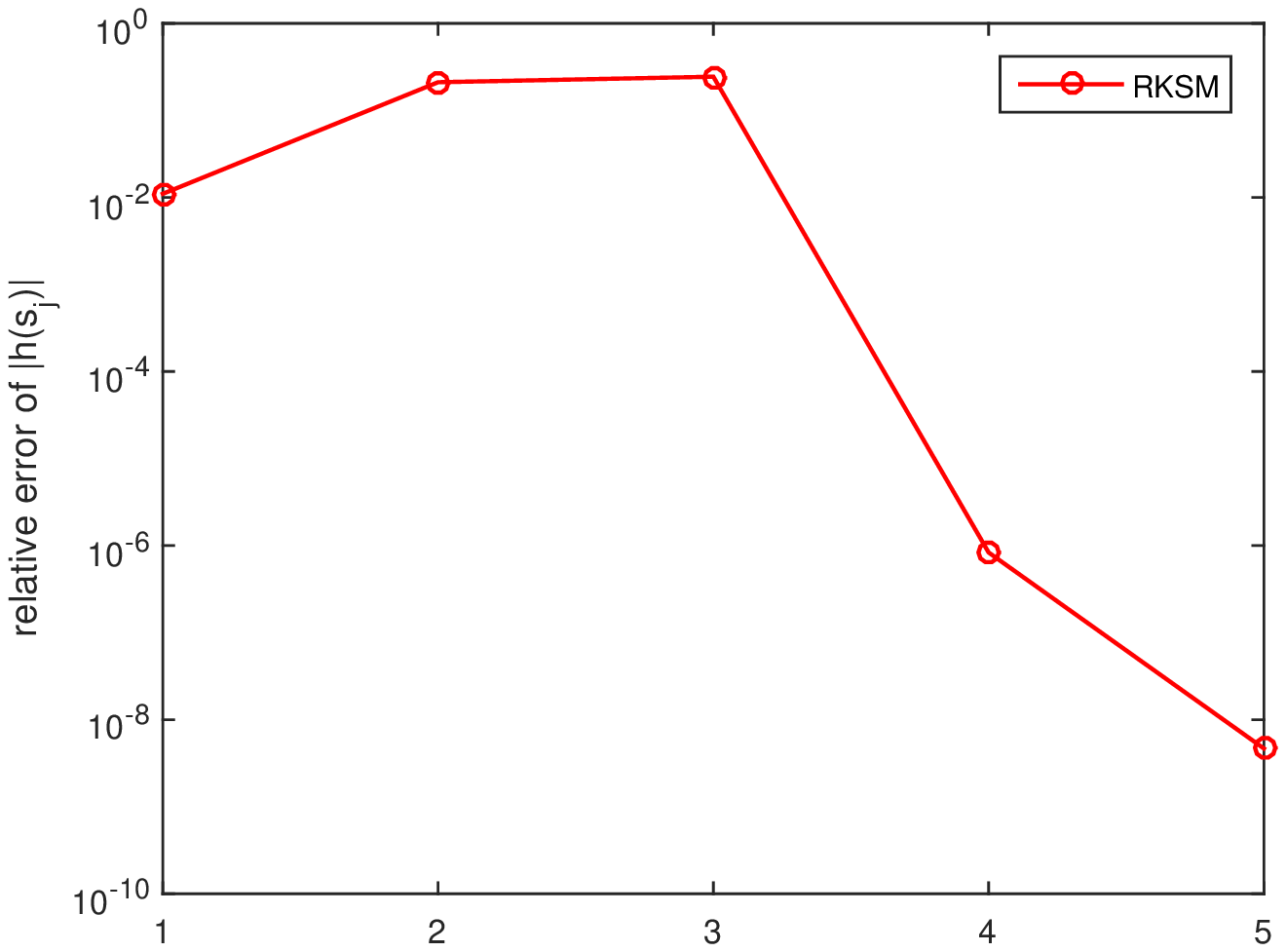}}

\footnotesize{In the left plot, the black asterisk  denotes  $(s_{*;1},h(s_{*;1}))$ and the red circle denotes the shifts $(s_j,h(s_j)) (j=3,4,5)$. As both $s_1$ and $s_2$ are outside the interval, they disappear.
The right  shows  the convergence of $|h(s_j)|$  to $0$.}
\end{figure}

To see the more detailed performance of RKSM, in Fig \ref{fig:kind2_1000_rc0.01}, we draw  $h(s)$ near $s_{*;1}$, and depict $|h(s_j)|$ when shifts ${s_{j(1:5)}}$ at Line 10 converges to $s_{*;1}$. One can see the fast convergence in Fig \ref{fig:kind2_1000_rc0.01}(b).

\subsection{Comparison with other algorithms}
In this subsection, we test several algorithms for  large-scale problems. 
Particularly, we compare RKSM  with six algorithms: \BN  ~\cite{zhya:2014a}, PsdLcp \cite{walz:2019}, \textsf{SDPT3} \cite{MR2894668,MR1976479},
\textsf{SeDuMi} \cite{MR1778433}, \textsf{cvx} \cite{gb08,cvx}
and LCPvA \cite{zhys:2015}.  \BN  ~is a Newton type method for  $h(s)=0$, and
LCPvA is a projection method \cite{zhys:2015} in a similar framework as RKSM.
Both \textsf{SDPT3} and \textsf{SeDuMi}, which have been included in the package \textsf{cvx}\footnote{{\tt http://cvxr.com/cvx/}.}, are  designed for general semidefinite-quadratic-linear programming. The procedure of transforming  an $\SOCLCP(\coneK^n,M,\bq)$ into a proper programming for \textsf{cvx} is given in \cite{MR1778433}, \cite[Section 4.6]{MR2894668} and \cite[Appendix]{zhys:2015}.
With $M=R^{\T}R$, the inputs for \textsf{SDPT3} and \textsf{SeDuMi} are $R$ and $\bq$, whereas  the inputs of \textsf{cvx} are $M$ and $\bq$.
In all the experiments, the CPU times of forming $M$ or $R$ are not counted in.



The settings of parameters in algorithms, except for RKSM, are the same as the ones in \cite{zhys:2015}; for RKSM, we set $\epsilon_1=10^{-7}$ and  $\epsilon_2=10^{-8}$ for the stopping criterion. Matrices of the following 
Example 2 and  Example 3 are the same as \cite[Section 7.2]{zhys:2015}. Due to version updating of some algorithms,   performances of relevant algorithms  change  slightly.

In our reported results,  ``--" means either the failure of obtaining a solution within given stopping criterion, or returning a solution with $\chi_{\rel}>1e-1$, where $\chi_{\rel}$ is defined by \eqref{eqn:relative_error}.
The label $d_M$ denotes the sparse densities of $M$; ``\# iter" represents the number of iterations, and 
{``CPU(s)" reports the output of Matlab \texttt{cputime}.
For our RKSM, ``\# iter" represents the dimension of the final subspace defined as
\begin{equation}\label{eqn:dim_iteration_relation}
\framebox{
\parbox{12cm}{
final dimension $=$ dimension of initial subspace $+$ $\ell$ $\times$ iteration number  $j$
}}
\end{equation}
In this subsection, our initial subspace is set as $\cK_{10,10}^{\ext}(J_nM,J_n\bq)$.
}

\vskip 3mm
\noindent\textbf{Example 2.}
This example is the same as  \cite[Table 7.3]{zhys:2015}. In \cite{zhys:2015}, it reports  average numbers of 4 results from 5 random examples. Here, we only present one result.  The  data matrix $M$ is again formed by $M=R^{\T}R,$ where $R$ is defined by (\ref{eq:sprandsym}) with $\texttt{n}=10000,
\texttt{desity}=0.0005$. The condition number of $M$ is about $(\frac{1}{\texttt{rc}})^2$ and can  range from $10^2$
to $10^5$. Results obtained from different \texttt{kind}s are listed in Table \ref{tab:large_sprandsym}, where \textsf{PsdLcp}  is excluded as  the storage is out of memory.

\def\vextras{\vphantom{\vrule height0.35cm width0.9pt depth0.15cm}}
\begin{table}[thb]
\renewcommand{\arraystretch}{1.0}
\caption{Numerical results of Example 2.}
\begin{center}\label{tab:large_sprandsym}
{\begin{small}  \tabcolsep0.05in
\begin{tabular}{|c|c|c|c|c|c|c|c|}
\hline
&& \multicolumn{3}{|c|}{\texttt{kind}=$1$} &\multicolumn{3}{|c|}{\texttt{kind}=$2$} \\ \cline{2-8}
&$\left(\frac{1}{\texttt{rc}}\right)^2(\approx \texttt{cond})$& $10^{2}$ &$10^{4}$ &$10^{5}$ & $10^{2}$ &$10^{4}$ &$10^{5}$  \vextras\\ \cline{2-8}
&$d_M$  & $7.3e\!-\!04$ & $7.4e\!-\!04$ & $7.5e\!-\!04$ & $2.6e\!-\!03$ & $2.6e\!-\!03$ & $2.6e\!-\!03$\\
 \hline
\multirow{6}{*}{\# iter}
&\small\BN             & $9/4$ & $9/4$ & $9/4$ & $13/2$ & $6/5$ & $4/9$ \\
&\small\textsf{SDPT3}  & $22$  & $19$  & $21$  & $15$   & $19$  & $21$  \\
&\small\textsf{Sedumi} & $16$  & $18$  & $16$  & $14$   & $19$  & $20$  \\
&\small\textsf{cvx}   & $19$	&$21$	&$22$&	$12$	&$15$	&$18$\\
&\small\LvA            & $30$  & $30$  & $30$  & $30$   & $50$  & $50$  \\
&\small\RKSM($\ell=1$)      & $ 22$&	$22$&	$22$&	$22$	&$23$	&$24$\\
&\small\RKSM($\ell=10$)     &$40$	&$40$	&$40$	&$40$	&$50$	&$50$\\
\hline
\multirow{6}{*}{CPU(s)}
&\small\BN             & $61.4$ & $60.7$ & $60.1$ & $112.9$ & $115.4$ & $173.2$ \\
&\small\textsf{SDPT3}  & $1.23$ & $1.03$ & $1.75$ & $57.8$  & $72.0$  & $80.4$  \\
&\small\textsf{Sedumi} & $1.10$ & $1.17$ & $1.11$ & $516.3$ & $718.3$ & $770.9$ \\
&\small\textsf{cvx} &$4.17$	&$2.90$	&$1.76$&	$2232.5$	&$2611.7$	&$3364.6$\\
&\small\LvA       & $0.14$&	$0.15$	&$0.15$&	$39.6$	&$61.9$&	$67.4$\\
&\small\RKSM($\ell=1$)   & $\textbf{0.10}$&	$\textbf{0.09}$&	$\textbf{0.09}$&	$\textbf{23.9}$&	$\textbf{30.3}$	&$\textbf{36.2}$\\
&\small\RKSM($\ell=10$) & $0.12$	&$0.11$	&$0.12$	&$38.9$	&$52.1$	&$52.4$\\

\hline
\multirow{6}{*}{$\chi_{\rel}$}
&\small\BN   & $6.4e\!-\!12$ & $7.3e\!-\!16$ & $1.0e\!-\!13$ & $6.1e\!-\!15$ & $2.1e\!-\!14$ & $1.2e\!-\!14$ \\
&\small\textsf{SDPT3}  & $1.9e\!-\!05$ & $9.5e\!-\!05$ & $1.2e\!-\!05$ & $3.8e\!-\!08$ & $2.3e\!-\!10$ & $8.7e\!-\!09$ \\
&\small\textsf{Sedumi} & $4.9e\!-\!07$ & $8.4e\!-\!07$ & $7.5e\!-\!07$ & $9.3e\!-\!08$ & $4.0e\!-\!09$ & $2.6e\!-\!10$ \\
&\small\textsf{cvx} &$4.8e\!-\!05$	&$3.3e\!-\!07$&	$1.3e\!-\!06$&	$3.5e\!-\!08$	&$1.6e\!-\!09$&	$9.6e\!-\!12$\\
&\small\LvA &$1.0e\!-\!12$	&$9.7e\!-\!13$&	$2.2e\!-\!13$&	$8.4e\!-\!10$&	$1.3e\!-\!11$&	$6.0e\!-\!08$\\
&\small\RKSM($\ell=1$)   &$9.6e\!-\!11$	&$1.1e\!-\!07$&	$9.8e\!-\!08$&	$1.4e\!-\!13$&	$1.0e\!-\!11$&	$5.7e\!-\!17$\\
&\small\RKSM($\ell=10$)   &$9.6e\!-\!11$	&$1.1e\!-\!07$&	$9.8e\!-\!08$&	$1.4e\!-\!13$&	$6.1e\!-\!17$&	$2.2e\!-\!13$\\

\hline
\end{tabular}
\end{small}
}\end{center}
\end{table}

From Table \ref{tab:large_sprandsym}, we observe that RKSM  converges  fastest. {For example, in 4 cases, RKSM($\ell=1$)  uses $22$ iterations;} since the dimension of the initial extended subspace is 20, this implies that  RKSM($\ell=1$) computes $s_{*;1}$ in one iteration, and  $s_{*;2}$ in another iteration.
The similar  discussion  applies to RKSM($\ell=10$)  {when the ``\# iter" number is $40$ (i.e., $40=20+10 \times 2$) by \eqref{eqn:dim_iteration_relation}.}


\vskip 3mm

\noindent\textbf{Example 3.}
This example is the same as  \cite[Table 7.4]{zhys:2015}. The tested  $M\succ 0$ are from Matrix Market, and we
do Cholesky decomposition $M=R^{\T}R$ and feed $R$ into \textsf{SDPT3} and \textsf{SeDuMi} as inputs. The results are showed in Table \ref{tal:matrix_market}.

\begin{table}[thb]
\renewcommand{\arraystretch}{1.2}
\caption{\small  Numerical results of Example 3.}
\begin{center}\label{tal:matrix_market}
{\begin{scriptsize}  \tabcolsep0.05in
\begin{tabular}{|c|c|c|c|c|c|c|c|c|c|c|}
\hline
&$M$& {\tiny \sc s1rmq4m1}&{ \tiny\sc s1rmt3m1}&{ \tiny\sc s2rmq4m1} &{\tiny \sc s2rmt3m1} &{\tiny \sc s3rmq4m1} &{\tiny \sc s3rmt3m1}&{\tiny \sc s3rmt3m3} &{\tiny \sc bcsstk17} &{\tiny \sc bcsstk18}  \\
\hline
&$n$              &$5489$        &$5489$ & $5489$&$5489$&$5489$&$5489$&$5357$&$10974$&$11948$\\
&$d_M$&$8.7e\!-\!03$&$7.2e\!-\!03$ &$8.7e\!-\!03$&$7.2e\!-\!03$&$8.7e\!-\!03$&$7.2e\!-\!03$&
$7.2e\!-\!03$&$3.6e\!-\!03$&$1.0e\!-\!03$\\
\hline
\multirow{7}{*}{\# iter}
&\scriptsize\BN               &  $25/1$&	$25/1$&	$21/1$	&$22/1$&	$18/2$	&$18/2$&	$18/2$	&$-$&	$-$   \\
&\scriptsize\textsf{PsdLcp}     &   $3$	&$3	$&$3$	&$3$	&$3$	&$3$&	$3$	&$-$&	$-$  \\
&\scriptsize\textsf{SDPT3} & $15$	&$16$	&$14$	&$13$	&$12$&	$18$	&$14$&	$19$&	$31$\\
&\scriptsize\textsf{Sedumi}& $20$	&$20$	&$17$	&$18$	&$16$	&$17$&	$17$	&$14$&	$20$\\
&\scriptsize\textsf{cvx}      &    $20$	&$20$&	$18$&	$18$	&$16$&	$15$&	$16$	&$11$&	$19$  \\
&\scriptsize\LvA               &$ 30$ &$30$ &$ 30$ &$30$ &$ 30$ &$ 30$ &$30$ &$20$ &$-$ \\
&\scriptsize\RKSM($\ell=1$)   &$24$&	$24$&	$24$	&$24$&	$24$&	$24$	&$24$	&$22$&	$34$\\
&\scriptsize\RKSM($\ell=10$)   &$50$	&$50$&	$50$& $50$	&$50$	&$50$	&$50$&	$40$&	$90$\\

\hline
\multirow{7}{*}{CPU(s)}
&\scriptsize\BN    &  $5.03$&	$4.10$&	$4.60$&	$3.88$	&$6.80$	&$6.32$&	$5.89$&	$-$	&$-$   \\
&\scriptsize\textsf{PsdLcp}    &      $ 213.7$&	$221.8$&	$244.6$&	$212.9$&	$213.7$&	$218.2$&	$196.7$&	$-$	&$-$  \\
&\scriptsize\textsf{SDPT3}& $14.38$	&$14.3$&	$12.6$&	$11.7$	&$11.1$	&$16.4$	&$326.0$	&$30.2$&	$243.2$\\
&\scriptsize\textsf{Sedumi}&$14.39$&	$14.9$	&$12.8$&	$13.2$&	$12.1$	&$12.9$	&$90.2$	&$19.9$	&$150.8$\\
&\scriptsize\textsf{cvx}  & $20.0$	&$16.5$&	$27.1$&	$16.3$&	$18.3$&	$8.20$&	$6.50$&	$11.2$&	$14.5$ \\
&\scriptsize\LvA   &$1.03$&	$0.61$	&$1.05$&	$0.57$&	$0.94$&	$0.58$&	$0.49$&	$1.42$&	$-$\\
&\scriptsize\RKSM($\ell=1$)&$\textbf{0.83}$&	$\textbf{0.51}$&	$\textbf{0.82}$&$	\textbf{0.51}$&	$\textbf{0.83}$&	$\textbf{0.53}$	&$\textbf{0.45}$&	$\textbf{0.72}$&	$\textbf{1.89}$\\
&\RKSM($\ell=10$) &$1.60$&	$0.97$&	$1.64$&	$0.95$	&$1.58$&	$0.99$&	$0.79$	&$1.41$&	$3.14$\\

\hline
\multirow{7}{*}{$\chi_{\rel}$}
&\scriptsize\BN  &$ 1.1e\!-\!13$	&$ 2.4e\!-\!13$&	$ 1.8e\!-\!12$&	$ 6.8e\!-\!12$&	$ 1.5e\!-\!16$&	$ 1.7e\!-\!17$
&	$ 1.1e\!-\!17$	&$-$&	$-$  \\
&\scriptsize\textsf{PsdLcp}	     & $ 4.7e\!-\!14$&	$ 7.8e\!-\!13$&	$ 8.1e\!-\!13$	&$ 1.8e\!-\!13$&	$ 1.6e\!-\!09$	&$ 2.3e\!-\!11$&	$ 2.2e\!-\!10$	&$-$&	$-$   \\

&\scriptsize\textsf{SDPT3}&$1.0e\!-\!06$&	$4.8e\!-\!07$&	$4.8e\!-\!08$&	$9.4e\!-\!07$&	$4.4e\!-\!07$&	$3.8e\!-\!10$	&$1.3e\!-\!07$&	$1.5e\!-\!16$	&$2.6e\!-\!14$\\
&\scriptsize\textsf{Sedumi}&$1.4e\!-\!08$&	$1.4e\!-\!08$	&$1.3e\!-\!09$&	$7.7e\!-\!09$	&$1.4e\!-\!07$	&$5.6e\!-\!08$	&$5.2e\!-\!08$	&$4.8e\!-\!17$	&$1.6e\!-\!15$\\
&\scriptsize\textsf{cvx} &$ 4.8e\!-\!05$	&$ 3.1e\!-\!05$&	$ 2.8e\!-\!06$&	$ 4.3e\!-\!06$&	$ 1.8e\!-\!06$
&	$ 1.5e\!-\!05$	&$ 3.2e\!-\!07$&	$ 3.1e\!-\!16$&	$ 6.0e\!-\!09$ \\

&\scriptsize\LvA &$ 3.0e\!-\!09$	&$ 4.0e\!-\!05$&	$ 2.7e\!-\!09$	&$ 1.0e\!-\!10$&	$ 6.7e\!-\!10$	&$ 1.2e\!-\!08$	&$ 1.4e\!-\!09$	&$ 5.3e\!-\!14$&	$-$\\

&\scriptsize\RKSM($\ell=1$)&$ 2.0e\!-\!10$	&$ 1.4e\!-\!11$	&$ 2.0e\!-\!10$&	$ 1.3e\!-\!11$	&$ 7.2e\!-\!11$&	
$ 1.3e\!-\!09$	&$5.1e\!-\!09$&	$ 2.4e\!-\!08$	&$ 3.7e\!-\!09$\\

&\scriptsize\RKSM($\ell=10$)&$ 6.1e\!-\!14$	&$ 4.1e\!-\!14$	&$ 1.6e\!-\!11$&	$ 2.6e\!-\!12$	&$ 3.5e\!-\!11$&
	$ 2.4e\!-\!10$	&$6.1e\!-\!09$&	$ 2.4e\!-\!08$	&$ 2.7e\!-\!16$\\

\hline
\end{tabular}
\end{scriptsize}
}\end{center}
\end{table}

We noticed that \textsf{PsdLcp} converges in only three iterations for most problems, but requires much more consuming time, mainly due to the full eigen-decompositions.
\textsf{LCPvA} obtains a lucky break on {\sc bcsstk17} but cannot deal with {\sc bcsstk18}, as reported in \cite{zhys:2015}. Due to a not good initial subspace, our RKSM($\ell=1$) spans a large subspace for {\sc bcsstk18} matrix $M$ with $\|M\|_1 \approx 5.12\times 10^{10},\tau \approx 3.61\times 10^{3}$ and $s_* \approx 3.88\times 10^{3}$. A close check of RKSM($\ell=1$) indicates that  the condition in line 7 is true, and thus  the new shifts $s_j=\|M\|_1/10^{j}$ are used in line 8. Until at $j=7$ (i.e., $\|M\|_1/10^{7}\approx 5.12\times 10^{3}$), the condition in Line 7 becomes false and Line 10 is executed.  Thus, the shifts are now near the target $s_* \approx 3.88\times 10^{3}$, and therefore, the  subspace formed becomes good enough for ensuring the convergence of  RKSM($\ell=1$).

By setting $\ell=10$ in RKSM for  {\sc bcsstk18}, the number of the \texttt{ldl} decompositions is reduced from 14 to 7, whereas the dimension of the final subspace increases from 34 to 90. In principal, it is possible to reduce the consuming CPU time as the callings of   \texttt{ldl} decompositions decrease.
However, due to the different efficiency of \texttt{backslash} and \texttt{ldl} for linear systems in Matlab,  for {\sc bcsstk18},  the reduced number of \texttt{ldl} decompositions is not enough to compensate the more expensive callings of   \texttt{ldl} decomposition than  \texttt{backslash}. The results are shown in Table \ref{tab:cpu_for_different_solvers}.  When the number of \texttt{ldl} calls for $\ell=1$ becomes sufficiently large,  the reduced number  (from using $\ell=10$) of \texttt{ldl}  may compensate the extra consuming time of \texttt{ldl}, and  we will see such an example in the {\sc flow} problem in  Example 4.


\begin{table}[htbp]
\centering
\caption { CPU  times of linear solvers for {\sc bcsstk18} in Example 3. \label{tab:cpu_for_different_solvers} }
\begin{small}
\begin{tabular}{|l|c|c|c|c|c|}
\hline
operation &   $M \backslash \bq$ & \texttt{ldl}$(M)$& \texttt{ldl}$(M), M \backslash \bq$ &
$\orth \left[\cK_{10}((J_nM)^{-1},(J_nM)^{-1}J_n\bq)\right]$\\ \hline
 CPU(s)    &  $0.054$ & $  0.164$& $ 0.180$ & $0.329$   \\ \hline
\end{tabular}
\end{small}
\end{table}



\vskip 3mm
\noindent\textbf{Example 4}.
We test larger problems in which the associated matrices $M$ are from the benchmark\footnote{{\tt https://sparse.tamu.edu/Oberwolfach}}\cite{10.1007/3-540-27909-1_11} for model order reduction. In the field of dynamical systems, there involves many \texttt{c}-stable $A$, whose eigenvalues are all on left hand  half plane. If \texttt{c}-stable $A$ is symmetric, then $A$ is negative positive definite. For those matrices, we simply take  $M=-A$ as test matrices.  The results are displayed in Table \ref{benchmark}.

\begin{table}[thb]
\renewcommand{\arraystretch}{1.2}
\caption{\small  Numerical results of Example 4.}
\begin{center}\label{benchmark}
{\begin{small}  \tabcolsep0.05in
\begin{tabular}{|c|c|c|c|c|c|c|c|c|c|}
\hline
&$M$& {\scriptsize\sc flow}&{\scriptsize \sc rail5177} &{\scriptsize \sc rail20209}  &{\scriptsize \sc gassensor}&{\scriptsize \sc chip} &{\scriptsize \sc t3dl} &{\scriptsize \sc t2dal} &{\scriptsize \sc t2dah} \\
\hline
 &$n$  &$ 9669$	&	$5177$&	$20209$	&	$66917$&	$20082$&	$20360$&	$4257$	&$11445$ \\
&$d_M$ &$7.2e\!-\!04$&		$1.3e\!-\!03$&	$3.4e\!-\!04$&	$3.8e\!-\!04$	&$7.0e\!-\!04$	&$1.2e\!-\!03$	&$8.7e\!-\!03$&	$1.3e\!-\!03$ \\
\hline

\multirow{6}{*}{\# iter}
&\scriptsize\textsf{PsdLcp} & $-$ & $2$   & $-$ & $-$ & $-$ & $-$ &$ 2 $  & $-$ \\
&\scriptsize\textsf{SDPT3}&$23$&	$-$&	$-$	&$-$&	$21$&	$-$	&$-$	&$-$\\
&\scriptsize\textsf{Sedumi}&$19$&	$14$	&$-$&	$-$	&$22$&	$18$&	$14$&	$-$\\
&\scriptsize\textsf{cvx}    & $18$  & $-$ & $-$ & $-$ & $-$ & $-$ & $-$ & $-$ \\
&\scriptsize\LvA  & $-$ & $30$  & $30$  & $-$ &$ 30$  & $-$ & $-$ & $-$ \\
&\scriptsize\RKSM($\ell=1$)& $32$&	$23$&	$23$	&$24$&	$23$&	$26$	&$24$	&$24$\\
&\scriptsize\RKSM($\ell=10$)&$60$&	$50$	&$50$	&$50$&	$50$&	$60$&	$50	$&$50$\\

\hline
\multirow{6}{*}{CPU(s)}
&\scriptsize\textsf{PsdLcp} & $-$  & $203.7$ & $-$  & $-$   & $-$  & $-$  & $109.0$ & $-$  \\
&\scriptsize\textsf{SDPT3}&$1031.6$&	$-$&	$-$	&$-$&	$1025.0$&	$-$&	$-$&	$-$\\
&\scriptsize\textsf{Sedumi}&$25.6$	&$146.2$&	$-$	&$-$&	$741.3$&	$901.2$&	$1.98$&	$-$\\
&\scriptsize\textsf{cvx}    & $6.30$  & $-$   & $-$  & $-$   & $-$  & $-$  & $-$   & $-$  \\
&\scriptsize\LvA & $-$	&$0.12$&	$0.60$&	$-$	&$8.41$	&$-$&	$-$	&$-$\\
&\scriptsize\RKSM($\ell=1$) & $0.82$&	$\textbf{0.09}$	&$\textbf{0.44}$&	$\textbf{85.0}$	&$\textbf{5.96}$&	$\scriptsize{18.3}$&	$\textbf{0.14}$&	$\textbf{0.57}$\\
&\scriptsize\RKSM($\ell=10$) & $\textbf{0.76}$&	$0.14$	&$0.85$&	$170.4$&	$13.8$	&$36.8$&	$0.19$	&$1.11$\\

\hline
\multirow{6}{*}{$\chi_{\rel}$}
&\scriptsize\textsf{PsdLcp} & $-$            & $ 3.6e\!-\!09$ & $-$            & $-$ & $-$            & $-$ & $ 7.1e\!-\!10$ & $-$ \\
&\scriptsize\textsf{SDPT3}&$ 1.7e\!-\!10$	&$-$	&$-$	&$-$	&$ 7.3e\!-\!05$	&$-$	&$-$	&$-$\\
&\scriptsize\textsf{Sedumi}&$ 4.5e\!-\!11$&	$ 7.8e\!-\!04$	&$-$&	$-$	&$ 8.6e\!-\!07$	&$ 1.9e\!-\!08$&	$ 8.6e\!-\!10$&	$-$\\
&\scriptsize\textsf{cvx}    & $ 6.1e\!-\!13$ & $-$            & $-$            & $-$ & $-$            & $-$ & $-$            & $-$ \\
&\scriptsize\LvA & $-$	&$ 1.4e\!-\!13$&	$ 7.4e\!-\!12$&	$-$&	$ 1.9e\!-\!12$&	$-$&	$-$	&$-$\\
&\scriptsize\RKSM($\ell=1$) &$6.5e\!-\!08$	&$ 7.6e\!-\!09$	&$ 8.2e\!-\!09$	&$ 5.8e\!-\!08$&	$ 9.7e\!-\!09$&	$3.6e\!-\!12$
&	$ 2.4e\!-\!16$	&$ 1.3e\!-\!16$\\

&\scriptsize\RKSM($\ell=10$) &$6.5e\!-\!08$	&$ 3.6e\!-\!09$	&$ 1.0e\!-\!10$	&$ 1.7e\!-\!09$&	$ 2.9e\!-\!09$
&	$ 6.0e\!-\!12$&	$ 1.3e\!-\!09$	&$ 1.0e\!-\!08$\\

\hline
\end{tabular}
\end{small}}
\end{center}
\end{table}

We noted that \textsf{BN} fails in all these problems, and \textsf{PsdLcp} solves two medium cases.
\textsf{cvx} successes in computing a solution within prescribed accuracy in only one example. By contrast, our RKSM works well for all these matrices. Also, for {\sc flow}, we observed that RKSM($\ell=10$) needs slightly less CPU time than  RKSM($\ell=1$) as the number of  \texttt{ldl} decompositions is reduced from  12 to 4.
%
%
%

\section{Conclusions}\label{section:conclusion}

Following the framework of \cite{zhys:2015}, in this paper, we proposed a new rational Krylov subspace method, RKSM($\ell$),  for solving large-scale symmetric and positive definite $\SOCLCP(\coneK^n,M,\bq)$. Through a transformation of $\SOCLCP(\coneK^n,M,\bq)$ into a zero-finding equation, we first connect it with the transfer functions in the model reduction. According to the moment match theory in model reduction, we observed in Theorem \ref{thm:h(s)-shifted} that the number of the  matched moments for $h(s)$ doubles when $M$ is symmetric. Thus, with a strategy of using multiple approximations, we propose RKSM, which improves the convergence and robustness over  \cite{zhys:2015} for the general $\SOCLCP(\coneK^n,M,\bq)$ with GUS property. Our numerical experiments demonstrate its efficiency and robustness.

\section*{Acknowledgments}
{The authors thank Dr. Ren-Cang Li at University of Texas at Arlington for discussions and comments on this paper.}

\bibliography{yd,strings,zhang-li}

\end{document}